\documentclass[11pt]{article}
\PassOptionsToPackage{obeyspaces}{url}
\usepackage{amsfonts, amsmath, amssymb}

\usepackage[hidelinks]{hyperref}
\usepackage{bookmark}
\bookmarksetup{open,numbered,depth=5} 

\usepackage{amsthm}

\usepackage{breakurl}
\usepackage{tikz}

\usepackage{etoolbox} 
\patchcmd{\thebibliography}{\leftmargin\labelwidth}{\leftmargin\labelwidth\addtolength\itemsep{-0.1\baselineskip}}{}{}

\usepackage{mathtools}

\newcommand*\samethanks[1][\value{footnote}]{\footnotemark[#1]}
\author{Boris Bukh\thanks{Department of Mathematical Sciences, Carnegie Mellon University, Pittsburgh, PA 15213, USA\@. Supported in part by U.S.\ taxpayers through NSF CAREER grant DMS-1555149. Email: {\tt bbukh@math.cmu.edu}, {\tt tchao2@andrew.cmu.edu}}
\and
Ting-Wei Chao\samethanks}

\title{Empty axis-parallel boxes}
\date{}

\usepackage[nameinlink]{cleveref}

\newtheorem{theorem}{Theorem}
\newtheorem{lemma}[theorem]{Lemma}

\newtheorem{definition}[theorem]{Definition}
\newtheorem{proposition}[theorem]{Proposition}
\newtheorem{claim}{Claim}
\crefname{claim}{Claim}{Claims}

\newcommand*{\eqdef}{\stackrel{\mbox{\normalfont\tiny def}}{=}}   
\newcommand*{\veps}{\varepsilon}                                 
\newcommand*{\R}{\mathbb{R}}                                     
\newcommand*{\Z}{\mathbb{Z}}                                     
\newcommand*{\E}{\mathbb{E}}                                     
\DeclarePairedDelimiter\abs{\lvert}{\rvert}                     
\DeclarePairedDelimiter\norm{\lVert}{\rVert}                     
\DeclareMathOperator{\vol}{vol}                                  
\DeclareMathOperator{\len}{len}                                  
\newcommand*{\T}{\mathcal{T}}                                    
\newcommand*{\Tor}{\mathbb{T}}                                   
\newcommand*{\A}{\mathcal{A}}                                    
\newcommand*{\D}{\mathcal{D}}                                    
\newcommand*{\Y}{\mathcal{Y}}                                    
\newcommand*{\bY}{\overline{\mathcal{Y}}}                        
\newcommand*{\cL}{\mathcal{L}}                                   

\begin{document}
\maketitle
\begin{abstract}
We show that, for every set of $n$ points in the $d$-dimensional unit cube, 
there is an empty axis-parallel box of volume at least~$\Omega(d/n)$ as $n\to\infty$ and $d$ is fixed.
In the opposite direction, we give a construction without an empty axis-parallel box of volume $O(d^2\log d/n)$. 
These improve on the previous best bounds of $\Omega(\log d/n)$ and $O(2^{7d}/n)$ respectively.
\end{abstract}

\section{Introduction}
\paragraph{Dispersion.} 
A \emph{box} is a Cartesian product of open intervals. Given a set $P\subset [0,1]^d$, we say that a 
box $B=(a_1,b_1)\times \dots \times (a_d,b_d)$ is \emph{empty} if $B\cap P=\emptyset$. Let $m(P)$
be the volume of the largest empty box contained in $[0,1]^d$.
Let $m_d(n)$ be the largest number such that every $n$-point set $P\subset [0,1]^d$ admits
an empty box of volume at least $m_d(n)$. Alternatively, $m_d(n)=\min m(P)$, where the minimum
is over all $n$-point sets $P\subset [0,1]^d$.

The quantity $m(P)$ is called the \emph{dispersion} of $P$. The motivation for estimating $m_d(n)$ came  independently in several subjects.
The earliest occurrence is probably in the work of Rote and Tichy~\cite{rote_tichy} who were
motivated by the relations to $\veps$-nets in discrete geometry on one hand, and
with the relations to discrepancy theory on the other. The dispersion
also arose in the problem of estimating rank-one tensors \cite{bddg,Krieg_Rudolf,nr_tensor} and in Marcinkiewicz-type
discretizations \cite{temlyakov}. In addition, in \cite{chen_dumitrescu}, lower bounds on dispersion were
used, via a compactness argument, to give constructions with large gap between covering and independence numbers for families of axis-parallel boxes.

The obvious bound $m_d(n)\geq 1/(n+1)$ was observed in several works, including \cite{dj60,bddg,rote_tichy}.
The first non-trivial lower bound of $m_d(n)\geq \tfrac{5}{4(n+5)}$ for $d\geq 2$ is due to Dumitrescu and Jiang
\cite{dj_amidst}. In \cite{dj60} they proved, for fixed $b$ and $d$, that $(n+1) m_d(n)\geq (b+1)m_d(b)-o(1)$,
which implies that the limit
\[
  c_d\eqdef \lim_{n\to\infty} n m_d(n)
\]
exists. Indeed, for each $b$, $\liminf (n+1) m_d(n)\geq (b+1)m_d(b)$, and therefore
$\liminf (n+1)m_d(n)\geq \limsup (n+1)m_d(n)$. Since $\liminf$ is always smaller than $\limsup$,
the limit exists.

The best lower bound on $m_d(n)$ for fixed $d$ is due to Aistleitner, Hinrichs and Rudolf \cite{ahr},
which is $c_d\geq \tfrac{1}{4}\log_2 d$. In the same paper they present a proof, due to Larcher,
that $c_d\leq 2^{7d+1}$. In this note we show that the correct dependence of $c_d$ on $d$ is neither logarithmic nor exponential, but polynomial.

\begin{theorem}\label{thm:crude}
The dispersion of $n$-point sets in $[0,1]^d$ satisfies
\begin{equation}
  m_d(n)\geq \frac{1}{n}\cdot \frac{2d}{e}\bigl(1-4dn^{-1/d}\bigr)\qquad\text{for all }d\text{ and all }n.\label{lowerbound:crude}
\end{equation}
\end{theorem}
All the logarithms in the rest of the paper are to base $e=2.718\dotsc$.
\begin{theorem}\label{thm:constr}
For every $d\geq 3$ and every $n\geq 1$, there is a set of at most $n$ points in $[0,1]^d$ for which the largest empty box
has volume at most $8000d^2\log d/n$.
\end{theorem}

For very large $n$, we have a slightly better lower bound.
\begin{theorem}\label{thm:sharp}
Let $R,T$ be positive real numbers that satisfy $T<R_0$ and $R_0-T<\log \frac{R_0}{T}$. Then
\[
  c_d\geq R_0\exp\bigl(-\tfrac{1}{2d}(R_0-T)\bigr).
\]
In particular, $c_d\geq \frac{2d}{e}(1+e^{-2d})$ for all $d$, and $c_2\geq 1.50476$.
\end{theorem}
This improves on the aforementioned bound of $c_2\geq 5/4$ by Dumitrescu--Jiang.
Very recently the upper bound of $c_2\leq 1.8945$ was proved by Kritzinger and Wiart \cite{kritzinger_wiart}.

\paragraph{Acknowledgements.} We thank Mario Ullrich and Daniel Rudolf for comments on the earlier version of this manuscript.

\section{Proofs of the lower bounds (\texorpdfstring{\Cref{thm:crude,thm:sharp}}{Theorems 1 and 3})}
\paragraph{Averaging argument.}
We first give a simple argument for \Cref{thm:crude}. We will then show how to modify
that argument to get \Cref{thm:sharp}. We start with the common part of the two arguments.\medskip

Let $R_0>0$ be a parameter to be chosen later subject to $R_0\leq n$, and set $\delta \eqdef \tfrac{1}{2}(R_0/n)^{1/d}$.
Let $f\colon [0,R_0]\to \R_+$ be some weight function. We postpone the actual choice of $f$ until later.
We adopt the convention that $f(R)=0$ if $R\geq R_0$.

Let $B$ be the cube of volume $R_0/n$ centered at the origin, i.e., $B\eqdef\bigl[-\delta,\delta\bigr]^d$.
Using $f$, we define a function on $\R^d$ by
\[
  F(x)\eqdef f(2^dr^dn)\qquad\text{ for }\norm{x}_{\infty}=r.
\]
Because $f$ vanishes outside $[0,R_0]$, the function $F$ vanishes outside~$B$.
Put $M\eqdef n\int_B F(x)\,dx$.
Note that $M=\int_{r=0}^{\delta} \bigl(2^dnr^{d-1}d\cdot f(2^dr^d n)\bigr)\,dr=\int_0^{R_0} f(R)\,dR$.
Because
\[
  \int_{t\in \R^d} \sum_{p\in P-t} F(p)\,dt=\sum_{p\in P}\int_{t\in \R^d} F(p-t)\,dt=\sum_{p\in P}\int_{x\in \R^d} F(x)\,dx=M,
\]
it follows that there exists $t\in [\delta,1-\delta\bigr]^d$ such that 
\begin{equation}\label{eq:weightbnd}
\sum_{p\in P-t} F(p)\leq M/(1-2\delta)^d,
\end{equation}
for otherwise $\int_{[\delta,1-\delta]^d} \sum_{p\in P-t} F(p)\,dt>M$.
It suffices to find a large box inside $B$ that is empty with respect to the set $P'\eqdef (P-t)\cap B$, for then we may obtain
an empty box of the same volume inside $[0,1]^d$ after translating by~$t$.


To find the empty box, we shave the sides off~$B$. Namely, for each point $p\in P'$
there is a coordinate of largest absolute value. If there is more than
one such coordinate, break the tie arbitrarily. Call this coordinate \emph{dominant}
for~$p$. Write the coordinates of $p\in P'$ as $p=(p_1,\dotsc,p_d)$. For each $i\in [d]$, put
\begin{align*}
  a_i&\eqdef\min \{-p_i : i\text{ is dominant for }p\in P'\text{ and }p_i\leq 0\},\\
  b_i&\eqdef\min \{\phantom{-}p_i : i\text{ is dominant for }p\in P'\text{ and }p_i\geq 0\}.
\end{align*}
Should the set in the definition of $a_i$ be empty, we put $a_i=\delta$.
Similarly, should the set in the definition of $b_i$ be empty, we put $b_i=\delta$.
The box 
\[
  B'\eqdef \prod_{i=1}^d (-a_i,b_i)
\]
is evidently disjoint from $P-t$ and is contained in~$B$. 

\begin{lemma}\label{lem}
The volume of $B'$ is at least $\frac{R_0}{n} \prod_{p\in P'} \sqrt{\frac{\norm{p}_{\infty}}{\delta}}$.
\end{lemma}
\begin{proof}
Fix any coordinate $i\in [d]$.

Suppose first that the two sets in the definitions of $a_i$ and $b_i$ are non-empty. 
Let $p,q\in P-t$ be the points such that $a_i=-p_i$ and $b_i=q_i$.
By the AM--GM inequality
\begin{equation}\label{eq:two}
  \frac{a_i+b_i}{2\delta}\geq \frac{\sqrt{a_ib_i}}{\delta}=\sqrt{\frac{\norm{p}_{\infty}}{\delta}}\cdot\sqrt{\frac{\norm{q}_{\infty}}{\delta}}.
\end{equation}

Suppose next that only one of the two sets in the definitions of $a_i$ and $b_i$ is non-empty.
Say $a_i=-p_i$ for some $p\in P-t$ and $b_i=\delta$ (the other case being symmetric). Then
by a similar application of the AM--GM inequality we obtain
\begin{equation}\label{eq:one}
  \frac{a_i+b_i}{2\delta}\geq \sqrt{\frac{\norm{p}_{\infty}}{\delta}}.
\end{equation}

By taking the product of \eqref{eq:two} and \eqref{eq:one} as appropriate over all $i\in [d]$, and noting that
every point has only one dominant coordinate,
we obtain
\[
  \vol B'=(2\delta)^d \cdot \prod_{i=1}^d \frac{a_i+b_i}{2\delta}\geq (2\delta)^d\prod_{p\in P'} \sqrt{\frac{\norm{p}_{\infty}}{\delta}}.\qedhere
\]
\end{proof}

\paragraph{Simple weight function (proof of \texorpdfstring{\Cref{thm:crude}}{Theorem 1}).}
The simplest choice of the constant $R_0$ and weight function $f$ is
\begin{align*}
  R_0&\eqdef 2d,\\
  f(R)&\eqdef \log \frac{R_0}{R}.
\end{align*}
The condition $R_0\leq n$ is satisfied unless $2d\leq n$, but in that case \Cref{thm:crude} holds vacuously.

With this choice of $R_0$ and $f$, we obtain $M=\int_0^{R_0} f(R)\,dR=R_0$ and $F(x)=d\log \frac{\delta}{\norm{x}_{\infty}}$ on $B$. 
So, from \Cref{lem} we obtain
\begin{align*}
  \vol B'&\geq \frac{R_0}{n} \exp\biggl(-\tfrac{1}{2}\sum_{p\in P'}\log \frac{\delta}{\norm{p}_{\infty}}\biggr)=
  \frac{R_0}{n} \exp\Bigl(-\tfrac{1}{2d}\sum_{p\in P'} F(p)\Bigr)\\
\intertext{which in view of \eqref{eq:weightbnd} is}
&\geq \frac{R_0}{n}\exp\Bigl(-\frac{1}{2d} M(1-2\delta)^{-d}\Bigr)=\frac{R_0}{n}\exp\Bigl(-\bigl(1-(R_0/n)^{1/d}\bigr)^{-d}\Bigr)\\
\intertext{which, since $(2d)^{1/d}\leq 2$, is}
&\geq \frac{R_0}{n}\exp\bigl(-(1-2n^{-1/d})^{-d}\bigr)\geq \frac{R_0}{n}\exp\bigl(-(1-2dn^{-1/d})^{-1}\bigr).\\
\intertext{Using $\exp(-(1-x)^{-1})=e^{-1}\cdot \exp(-x-x^2-\dotsb)\geq e^{-1}\cdot (1-x-x^2-\dotsb)\geq e^{-1}(1-2x)$ for $x\in[0, 1/2]$, we may deduce that}
\vol B' &\geq \frac{1}{n}\cdot\frac{2d}{e}\bigl(1-4dn^{-1/d}\bigr).
\end{align*}

\paragraph{Better weight function (proof of \texorpdfstring{\Cref{thm:sharp}}{Theorem 3}).}
Let $T$ and $R_0$ be as in the statement of \Cref{thm:sharp}.
Since the aim is to prove a bound on $c_d$, we may assume that $n$ is sufficiently large.
Define
\begin{equation}
  f(R)\eqdef 
    \begin{cases}
      \log \frac{R_0}{T}&\text{if }R\leq T,\\
      \log \frac{R_0}{R}&\text{if } T<R\leq R_0.\\       
    \end{cases}
\end{equation}
It is readily computed that $M=(\int_0^T+\int_T^{R_0})f(R)\,dR=R_0-T$. Since $(1-2\delta)^d\to 1$, it follows that 
$M/(1-2\delta)^d<\log \frac{R_0}{T}$, for large enough $n$. 
Because of \eqref{eq:weightbnd}, this implies that
$\sum_{x\in P'} F(x)<\log \frac{R_0}{T}$, and hence
for no point $x\in P'$ does it hold that $R\leq T$, where $R=2^dn\norm{x}_{\infty}^d$.
So, $F(x)=d\log \frac{\delta}{\norm{x}_{\infty}}$ for all $x\in P'\cap B$.
So, we may proceed as before to obtain
\begin{align*}
  \vol B'&\geq \frac{R_0}{n} \exp\Bigl(-\tfrac{1}{2d}\sum_{p\in P'} F(p)\Bigr)\geq \frac{R_0}{n} \exp\bigl( -\tfrac{1}{2d}M(1-2\delta)^{-d}\bigr).
\end{align*}
Taking the limit $n\to\infty$, the bounds on $c_d$ follows.\smallskip

The bound $c_d\geq \frac{2d}{e}(1+e^{-2d})$ is obtained by choosing $R_0=2d$ and $T=R_0\exp(-R_0)$. The bound
$c_2\geq 1.50476$ is obtained by choosing $R_0=3.69513$ and $T=0.101622$.

\section{Proof of the upper bound (\texorpdfstring{\Cref{thm:constr}}{Theorem 2})}
\paragraph{Construction outline.}
Our construction is a modification of the Hilton--Hammerseley construction. As in the
Halton--Hammerseley construction, we will select primes $p_1,\dotsc,p_d$, each of which is associated to respective coordinate direction.
As in the analysis of Halton--Hammerseley construction, we will be interested in \emph{canonical boxes}, which are the boxes\footnote{Here and elsewhere in this section we work with half-open boxes.
Since every half-open box contains an open box of the same volume, this does not impair the strength of our constructions, but
doing so will be technically advantageous.} of the form
\[
  B=\prod_{i=1}^d \left[\frac{a_i}{p_i^{k_i}},\frac{a_i+1}{p_i^{k_i}}\right).
\]
for some integers $0\leq a_i<p_i^{k_i},i=1,2,\dotsc,d$.

For a prime $p$ and a nonnegative integer $x$, consider the base-$p$ expansion of the number $x$, say
$x=\nobreak x_0+x_1p+\dotsb+x_{\ell}p^{\ell}$. Put $r_p(x)\eqdef x_0p^{-1}+x_1p^{-2}+\dotsb+x_{\ell}p^{-\ell-1}$; note that
$r_p(x)$ is the number in $[0,1)$ obtained by reversing the base-$p$ digits of~$x$.
Define the function $r\colon \Z_{\geq 0}\to [0,1]^d$ by $r(x)\eqdef \bigl(r_{p_1}(x),\dotsc,r_{p_d}(x)\bigr)$.

Our construction is broken into two stages. The set that we construct in the first stage
is an $r$-image of a certain subset of $\Z_{\geq 0}$. (Note that the usual Halton--Hammerseley construction is the $r$-image of an interval of length~$n$.)
This set has $O(n d \log d)$ elements and intersects \emph{almost} all the boxes of volume about~$1/n$. 
In the second stage of the construction, we show that $d+1$ suitably chosen translates 
of the first set meet all the boxes of volume~$1/n$.

\paragraph{First stage.} To simplify the proof, we will discretize the boxes we work with. We will do so by shrinking them slightly,
so that $i$'th coordinates have terminating base-$p_i$ expansions. 

With hindsight we choose $p_i$ to be the $(d+i)$'th smallest prime, for each $i=1,2,\dots,d$. 
Put $\gamma\eqdef p_1p_2\dots p_d$, and let $n$ be an arbitrary integer divisible by $2\gamma^{11}$.

\begin{definition}
We say that a box $\beta$ is a \emph{good box} if   it is of the form
\begin{equation}\label{eq:goodboxform}
\beta=\prod_{i=1}^d \left[\frac{a_i}{p_i^{k_i}}+\frac{b_i}{p_i^{k_i+3}},\frac{a_i}{p_i^{k_i}}+\frac{c_i}{p_i^{k_i+3}}\right),
\end{equation}
for some integers $0\leq b_i<c_i\leq p_i^3$ and $k_i\in\Z_{\geq 0}$ for $i=1,2,\dots,d$, and whose volume is $1/4n\leq\nobreak \vol(\beta)\leq 1/n$. 
Let $B=\prod_i \left[a_i/p_i^{k_i},(a_i+1)/p_i^{k_i}\right)$.
We call $(B,\beta)$ a \emph{good pair}.
\end{definition}
Since $i$'th coordinate dimension of $B$ is at most $p_i^3$ times larger than that of $\beta$, it follows that $\vol(B)\leq \gamma^3/n$.
In other words, a (discretized) box $\beta$ is good if it is contained in a canonical box $B$ that is not much larger than $\beta$.
Note that, since a good $\beta$ can sometimes be written in the form \eqref{eq:goodboxform} in more than one way,
the choice of $B$ in the definition of a good pair is, in general, not unique.\medskip

Our aim in this stage of construction is to find a set $P$ that meets every good box. In the next stage
we will superimpose several copies of $P$ to create a set that meets every large box. It is precisely
because the family of good boxes is richer than the family of canonical boxes
that we lose less in the second stage than if we used the Halton--Hammerseley construction.
\medskip

Suppose $B$ is a canonical box. Write it as $B=\prod_i\left[a_i/p_i^{k_i},(a_i+1)/p_i^{k_i}\right)$, and consider $r^{-1}(B)$. The set $r^{-1}(B)$ consists
of the solutions to the system
\begin{align*}
x&\equiv a_1'\pmod{p_1^{k_1}},\\
x&\equiv a_2'\pmod{p_2^{k_2}},\\
 &\setbox0\hbox{$\equiv$}\mathrel{\makebox[\wd0]{\vdots}}\\
x&\equiv a_d'\pmod{p_d^{k_d}},
\end{align*}
where $a_i'\eqdef r_{p_i}(a_i)p_i^{k_i}$, i.e., $a_i'$ is the integer obtained from $a_i$ by reversing its base-$p_i$ expansion. 

By the Chinese Remainder theorem, the set $r^{-1}(B)$ is an infinite arithmetic progression with step $D(B)\eqdef p_1^{k_1}p_2^{k_2}\dotsb p_d^{k_d}=1/\vol(B)$.
Let $A(B)$ be the least element of $r^{-1}(B)$, so that 
\[r^{-1}(B)=\lbrace A(B)+kD(B): k\in\mathbb{Z}_{\geq 0}\rbrace.\]
Given a good pair $(B,\beta)$, define \[L_B(\beta)\eqdef\{ k \in\mathbb{Z}_+ : r\bigl(A(B)+kD(B)\bigr) \in \beta\}.\]

\begin{claim}\label{claim:L}
The set $\cL\eqdef \lbrace L_B(\beta):(B,\beta)\mbox{ is a good pair}\rbrace$ is of size at most $\gamma^{12}$.
\end{claim}

\begin{proof}
Let $(B,\beta)$ be a good pair. Write $B$ and $\beta$ in the form
\[B=\prod_{i=1}^d \left[\frac{a_i}{p_i^{k_i}},\frac{a_i+1}{p_i^{k_i}}\right),\qquad \beta=\prod_{i=1}^d \left[\frac{a_i}{p_i^{k_i}}+\frac{b_i}{p_i^{k_i+3}},\frac{a_i}{p_i^{k_i}}+\frac{c_i}{p_i^{k_i+3}}\right).\]

We know that $r\bigl(A(B)+kD(B)\bigr)\in \beta$ is equivalent to
\begin{align*}
A(B)+kD(B)&\in a_1'+p_1^{k_1}J_1\pmod{p_1^{k_1+3}},\\
A(B)+kD(B)&\in a_2'+p_2^{k_2}J_2\pmod{p_2^{k_2+3}},\\
&\setbox0\hbox{$\in$}\mathrel{\makebox[\wd0]{\vdots}}\\
A(B)+kD(B)&\in a_d'+p_d^{k_d}J_d\pmod{p_d^{k_d+3}},
\end{align*}
where the sets $J_i$ consist of base-$p_i$ reversals of the numbers in the interval $[b_i,c_i)$ (which are $3$-digit long in base $p_i$).

On the other hand, we know that
\begin{align*}
A(B)+kD(B)&\equiv a_1'+(\alpha_1+k\delta_1)p_1^{k_1}\pmod{p_1^{k_1+3}},\\
A(B)+kD(B)&\equiv a_2'+(\alpha_2+k\delta_2)p_2^{k_2}\pmod{p_2^{k_2+3}},\\
&\setbox0\hbox{$\equiv$}\mathrel{\makebox[\wd0]{\vdots}}\\
A(B)+kD(B)&\equiv a_d'+(\alpha_d+k\delta_d)p_d^{k_d}\pmod{p_d^{k_d+3}}
\end{align*}
for some $\alpha_i,\delta_i\in\mathbb{Z}/p_i^3\mathbb{Z},i=1,2,\dots,d$. There are at most $\gamma^6$ different choices for $(\alpha_i,\delta_i)_{i=1}^d$.
Also, there are at most $\gamma^6$ different choices for $(b_i,c_i)_{i=1}^d$ satisfying $0\leq b_i<c_i\leq p_i^3$. Since $L_B(\beta)$ is determined by $(\alpha_i,\delta_i,b_i,c_i)_{i=1}^d$, the claim is true.
\end{proof}

To each canonical box $B$ of volume between $1/4 n$ and $\gamma^3/n$ we assign
a \emph{type}, so that boxes of the same type behave similarly.
Formally, let $\A(B)$ be the unique multiple of $n/\gamma^4$ satisfying 
$0\leq A(B)-\A(B)<n/\gamma^4$.
Similarly,
let $\D(B)$ be the unique multiple of $n/\gamma^{11}$ satisfying \linebreak $0\leq D(B)- \D(B)< n/\gamma^{11}$.
The type of $B$ is then the pair $\T(B)\eqdef \bigl(\A(B),\D(B)\bigr)$. 

Note that, from $1/4n\leq \vol(B)\leq \gamma^3/n$ and $D(B)=1/\vol(B)$ it follows that
\begin{equation}\label{eq:dbound}
  n/\gamma^3-n/\gamma^{11}<\D(B)\leq 4n.
\end{equation}

\begin{claim}\label{claim:type}
The number of types is at most $\gamma^{16}$.
\end{claim}
\begin{proof}
Since $A(B)<D(B)\leq 4 n$, the number of types is at most
$(\frac{4 n}{n/\gamma^4})(\frac{4 n}{n/\gamma^{11}})=16\gamma^{15}\leq \gamma^{16}$.
\end{proof}

For a type $\T=(\A,\D)$, let $\Y(\T)\eqdef \lbrace \A+k\D : k\in\mathbb{Z}_{\geq 0}\rbrace$ be the arithmetic progression generated by 
$\A$ and $\D$. Note that if $\T=\T(B)$, then $\Y(\T)$ is an approximation to $r^{-1}(B)$.
In particular, $\Y(\T)$ and $r^{-1}(B)$ intersect any long interval that is not too far from the origin in approximately
the same way.

For integers $a,b$, denote by $[a,b)$ the integer interval consisting of integers $x$ satisfying $a\leq x<b$.
Our construction will be a union of intervals of length $n/\gamma^3$ whose left endpoints are in $[0,n\gamma^4)$.\smallskip

We first estimate the difference between respective terms in $\Y(\T)$ and $r^{-1}(B)$ inside $[0,n\gamma^4)$.
\begin{claim}\label{claim:shrink}
Suppose $\T(B)=\bigl(\A(B),\D(B)\bigr)$.
Then for any integer $x\in [0,n\gamma^4)$ and any integer $k$,
$\A(B)+k\D(B)\in [x,x+ n/2\gamma^3)$ implies $A(B)+kD(B)\in [x,x+ n/\gamma^3)$.
\end{claim}

\begin{proof}
For such $k$, since $\mathcal{A}(B)+k\mathcal{D}(B)< n\gamma^4+ n/\gamma^3$, from \eqref{eq:dbound} we deduce that \[k< \frac{n\gamma^4+ n/\gamma^3}{n/\gamma^3-n/\gamma^{11}}\leq 2\gamma^7.\] 
In view of $k\geq 0$, this implies that
\[0\leq (A(B)+kD(B))-(\mathcal{A}(B)+k\mathcal{D}(B))\leq \frac{n}{\gamma^4}+2\gamma^7\cdot \frac{n}{\gamma^{11}}= \frac{3n}{\gamma^4},\]
and hence $A(B)+kD(B)\in[x,x+n/2\gamma^3+3n/\gamma^4)\subseteq [x,x+n/\gamma^3)$.
\end{proof}

For a type $\T$ and $L\in\cL$ that satisfy $\T=\T(B)$ and $L=L_B(\beta)$ for some good pair $(B,\beta)$, define
\[\Y_{\T}(L)\eqdef\lbrace \A+k\D:k\in L\rbrace.\]
With this definition,
$\Y_{\T}(L)$ 
is the approximation to $r^{-1}(\beta)$ induced by the approximation $\Y(B)$
to~$r^{-1}(B)$.

\begin{claim}
The set $\bY_{\T}(L)=\Y_{\T}(L)\cap [0,n\gamma^4)$ is of size at least $\gamma^4/16+1$.
\end{claim}

\begin{proof}
Let $(B,\beta)$ be a good pair such that $\T=\T(B)$ and $L=L_B(\beta)$.
The set $L_B(\beta)$ is $\gamma^3$-periodic, i.e., $k\in L_B(\beta)$ implies $k+\gamma^3\in L_B(\beta)$. 
The intersection of any interval of length $\gamma^3$ with $L_B(\beta)$ is of size exactly $\gamma^3\frac{\vol(\beta)}{\vol(B)}$.
Since the preimage of $[0,n\gamma^4)$ under the map $k\mapsto \A+k\D$ contains 
\[
  \left\lfloor \frac{n\gamma^4-\A}{\gamma^3\D}\right\rfloor\geq \frac{n\gamma}{\D}-2\geq n\gamma \vol(B)-2\geq \frac{1}{2}n\gamma \vol(B)
\]
non-overlapping intervals of length $\gamma^3$, the size of $\bY_{\T}(L)$ is at least 
\[\frac{1}{2}n\gamma \vol(B)\cdot \gamma^3\frac{\vol(\beta)}{\vol(B)}=\frac{1}{2}n  \gamma^4 \vol(\beta)\geq \gamma^4/16+1.\qedhere\]
\end{proof}

\begin{claim}\label{claim:probability}
Let $x$ be chosen uniformly from $[0,n\gamma^4)$. Then
$\Pr\bigl[\,\bY_{\T}(L)\cap [x,x+n/2\gamma^3) \neq \emptyset \bigr]\geq 1/32\gamma^3$.
\end{claim}
\begin{proof}
Let $y\in \bY_{\T}(L)$ be arbitrary. If $y\notin [0,n/2\gamma^3)$, then $\Pr[y\in [x,x+n/2\gamma^3)]=1/2\gamma^7$. Since
$\D>n/\gamma^3-n/\gamma^{11}\geq n/2\gamma^3$, the set $\bY_{\T}(L)$ contains at most one element in the interval $[0,n/2\gamma^3)$.
Hence \[\E\bigl[ \abs[\big]{\bY_{\T}(L)\cap [x,x+n/2\gamma^3)} \bigr]\geq 1/32\gamma^3.\]
Since elements of $\bY_{\T}(L)$ are at least $\D$ apart, $\abs{\bY_{\T}(L)\cap [x,x+n/2\gamma^3)}\in \{0,1\}$ for all $x$.
Therefore, \[\Pr\bigl[\,\bY_{\T}(L)\cap [x,x+n/2\gamma^3) \neq \emptyset \bigr]=\E\bigl[ \abs[\big]{\bY_{\T}(L)\cap [x,x+n/2\gamma^3)} \bigr]\geq 1/32\gamma^3.\qedhere\]
\end{proof}

Sample $900 \gamma^3 \log \gamma$ elements uniformly at random from $[0,n\gamma^4)$, independently from one another.
Let $X$ be the resulting set. Then by the preceding claim
\begin{align*}
  \Pr\bigl[\,\bY_{\T}(L)\cap \bigl(X+[0,n/2\gamma^3)\bigr)=\emptyset\bigr]\leq (1-1/32\gamma^3)^{900 \gamma^3 \log \gamma}<\gamma^{-28}.
\end{align*}
From \Cref{claim:L,claim:type} and the union bound it then follows that there exists a choice of $X$
such that $\bY_{\T}(L)\cap \bigl(X+[0,n/2\gamma^3)\bigr)$ is non-empty whenever $\T=\T(B)$, $L=L_B(\beta)$ and $(B,\beta)$ is a good pair.
In other words, for every $(B,\beta)$ there exist $x\in X$ and an integer $k\in L_B(\beta)$ such that $\A(B)+k\D(B)\in [x,x+n/2\gamma^3)$.
By \Cref{claim:shrink} this implies that $A(B)+kD(B)\in [x,x+n/\gamma^3)$ for the same $x$ and~$k$, whereas the definition
of $L_B(\beta)$ implies that $r\bigl(A(B)+kD(B)\bigr)\in \beta$.
Because this holds for every good pair $(B,\beta)$,
the set $P\eqdef r\bigl(X+[0,n/\gamma^3) \bigr)$ meets every good box.

Note that
$\abs{P}\leq \abs{X}\cdot \frac{n}{\gamma^3}\leq 900 \log \gamma\cdot n\leq 3000 d\log d \cdot n$ (since $\log \gamma\leq d\log p_d\leq 3d\log d$).

\paragraph{Second stage.} So far we have worked with boxes whose coordinates are rational numbers with denominators
of the form $p_i^{k_i}$. Given an arbitrary box, we shall shrink it down to a box of such form. We begin by describing this
process.

A \emph{$p$-interval} is an interval of the form $[a/p^k,b/p^k)$ for some integers $0\leq a<b<p^k$.
A \emph{canonical $p$-interval} is an interval of the form $[a/p^k,(a+1)/p^k)$ with $0\leq a<p^k$.
Note that canonical boxes are precisely the boxes that are Cartesian products of canonical intervals in appropriate bases.
A $p$-interval $[a/p^k,b/p^k)$ is \emph{well-shrunk} if $b-a<p^2$. 

\begin{claim} Every interval $[s,u)$
contains a well-shrunk $p$-interval of length at least $(1-2/p)\len [s,u)$.
\end{claim}
\begin{proof}
Let $k$ be the smallest integer satisfying $\len [s,u)\geq p^{-k}$. Let $I$ be the largest interval
of the form $I=[a/p^{k+1},b/p^{k+1})$ contained in $[s,u)$. Then 
$\len I
\geq u-s-2p^{-k-1} \geq (1-2/p)(u-s)$, and
$b-a=p^{k+1}\len I\leq p^{k+1}\len [s,u)<p^2$.
\end{proof}

Call an interval $[s,u)$ \emph{$p$-bad} if it contains a rational number with denominator $p^{k+1}$, where $\len [s,u)<2p^{-k-2}$ and $k\in \Z_{\geq 0}$.
\begin{claim}\label{claim:bad}
A box $\alpha=\prod_i[s_i,u_i)\subset [0,1]^d$ of volume $1/n$ fails to contain a good box only if, for some $i\in[d]$,
the interval $[s_i,u_i)$ is $p_i$-bad.
\end{claim}
\begin{proof}
For each $i\in [d]$, let $[s_i',u_i')$ be a well-shrunk $p_i$-interval
contained in $[s_i,u_i)$ as above. Let $\beta\eqdef \prod_i [s_i',u_i')$. Note that $\vol(\beta)\geq \vol(\alpha)\prod_i (1-2/p_i)\geq 1/4n$.

Let $B=\prod_i [a_i/p_i^{k_i},(a_i+1)/p_i^{k_i})$ be the smallest canonical box containing $\beta$. 
Since the $p$\nobreakdash-interval $[s_i',u_i')$ is contained in $[a_i/p_i^{k_i},(a_i+1)/p_i^{k_i})$, we may write
it in the form \[[s_i',u_i')=[a_i/p_i^{k_i}+b_i/p_i^{\ell_i},a_i/p_i^{k_i}+c_i/p_i^{\ell_i})\] for
some integers $0\leq b_i<c_i<p_i^{k_i-\ell_i}$. Since $[s_i',u_i')$ is well-shrunk, $c_i-b_i<p^2$.

If $(B,\beta)$ is not a good pair, there exists $i\in [d]$, such that $\ell_i\geq k_i+4$. Fix such an $i$.
By the minimality of $B$, the interval $[s_i',u_i')$ contains a rational number with denominator
$p_i^{k_i+1}$. Since $[s_i,u_i)$ contains $[s_i',u_i')$, this rational number is also contained
in $[s_i,u_i)$. As $\len [s_i,u_i)\leq (c_i-b_i+2)p^{-\ell_i}<(p^2+2)p^{-k_i-4}\leq 2p^{-k_i-2}$, the interval $[s_i,u_i)$ is $p_i$-bad.
\end{proof}
\begin{claim}\label{claim:translates}
Let $\Delta=1/p(p-1)$. Suppose $[s,u)\subset [1/p,1]$ is an arbitrary interval. Then at most one of its translates $[s,u),[s,u)-2\Delta,\dotsc,[s,u)-2d\Delta$ is $p$-bad.
\end{claim}
\begin{proof}
Suppose that, for some $r$, the interval $[s,u)-2r\Delta$ contains rational number $a/p^{k+1}$ and is of length $\len [s,u)<2p^{-k-2}$.
Then the interval $[s_r,u_r)\eqdef [s,u)- 2r\Delta - a/p^{k+1}$ contains $0$ and is also of length $\len I<2p^{-k-2}$.
Hence, $u_r< 2p^{-k-2}$, and so $(k+2)$'nd digit in the base-$p$ of $u_r$ is either $0$ or~$1$.
Note that it is the same as the $(k+2)$'nd digit of $u-2r\Delta$.

Since the base-$p$ expansion of $\Delta$ is $0.01111\dotsb$ and $2d+1<p$, 
for at most one of the numbers $u,u-2\Delta,\dotsc,u-2d\Delta$ is the $(k+2)$'nd digits
equal to $0$ or~$1$.
Hence, at most one of the intervals $[s,u),[s,u)-2\Delta,\dotsc,[s,u)-2d\Delta$ contains a rational number with denominator $p^{k+1}$.
\end{proof}

Let $P$ be the set constructed in the first stage. Let $v\in [0,1]^d$ be the vector whose $i$'th coordinate
is $v_i=1/p_i(p_i-1)$. Let $P'\eqdef \bigcup_{r=0}^d (P+2rv)$. We claim that $P'$ meets every subbox
of $\prod_i [1/p_i,1]$ of volume~$1/n$.

Indeed, suppose $\alpha=\prod_i [s_i,u_i) \subset \prod_i [1/p_i,1]$ is an arbitrary box of volume~$1/n$.
Then by the preceding claim, there exists $r\in \{0,1,\dotsc,d\}$ such that for no $i\in[d]$ is the interval $[s_i,u_i)-2r\Delta$ 
$p$-bad. \Cref{claim:bad} tells us that the box $\alpha-2r\Delta$ contains a good box. Since
the set $P$ meets all good boxes, it follows that the $P+2r\Delta$ meets $\alpha$.
As $P+2r\Delta\subset P'$, the set $P'$ indeed meets~$\alpha$.

Finally, we scale the box $\prod_i [1/p_i,1]$ onto $[0,1]^d$. This way, we turn the set $P'$
into a set that meets every subbox of $[0,1]^d$ of volume $\frac{1}{n}\prod p_i/(p_i-1)\leq 2/n$.
This set has size $\abs{P'}\leq (d+1)\cdot 3000d\log d\cdot n$.\medskip

This construction shows that $m_d(\lfloor 3000d(d+1)\log d\cdot n\rfloor)\leq 2/n$ for all $n$ that are divisible by $2\gamma^{11}$.
Since the limit $c_d=\lim_{n\to\infty} nm_d(n)$ exists, it then follows that
$c_d\leq 6000d(d+1)\log d$, which, by the Dumitrescu--Jiang inequality mentioned in the introduction,
implies that $m_d(b)\leq  \frac{1}{b+1}\cdot 6000d(d+1)\log d$ for all $b$.
Because $6000d(d+1)\log d\leq 8000d^2\log d$, the proof is complete.



\section{Problems and remarks}
\begin{itemize}
\item Because of the $n^{-1/d}$ term, the bound in \Cref{thm:crude} is weak when the number
of points $n$ is small compared to the dimension~$d$. It is likely possible to
replace the term $n^{-1/d}$ with $O_d(n^{-1})$ by using a more sophisticated averaging argument. 
In our argument we considered an average of translates of a
function supported on a fixed box~$B$. The error term $n^{-1/d}$ is due to
the points near the boundary of $[0,1]^d$ receiving less weight than the rest.
One can remedy this by using, in addition to the translates of $B$,
also elongated boxes of volume $\vol(B)$ to add weight in the regions
near the boundary of $[0,1]^d$. In this paper, we decided to sacrifice 
the slightly stronger bound for a simpler proof.

For best constructions of low-dispersion sets that are good when $n$ is small compared to $d$, see \cite{ullrich_vybiral,litvak} (improving
on the earlier bounds in \cite{sosnovec}). These constructions are probabilistic. For an explicit construction, which has larger dispersion, see \cite{krieg}.

\item The low-dispersion sets are used in \cite{bddg}, \cite{Krieg_Rudolf}, and \cite[Theorem 11]{nr_tensor} to give algorithms
to approximate certain one-dimensional tensors. Because of that, it is useful to
derandomize the construction in \Cref{thm:constr}. The following is a way to do so. It gives an algorithm that computes a set from \Cref{thm:constr}
in $d^{O(d)}+O(dn)$ arithmetic operations.

The algorithm is broken into two steps. The first step is a pre-processing step, which depends solely on~$d$.
The second step takes the output of the first step and $n$ and quickly produces the $n$-point low-dispersion point set in $[0,1]^d$.

For the pre-processing step, we need some definitions.

\begin{definition}
For any tuple $(\alpha_i,\delta_i,b_i,c_i)_{i=1}^d$ with $\alpha_i,\delta_i,b_i,c_i\in\mathbb{Z}/p_i^3\mathbb{Z},i=1,2,\ldots,d$, consider the system of $d$ equations in unknown $k$
\begin{align*}
\alpha_1+k\delta_1&\in J_1 \pmod{p_1^{3}},\\
\alpha_2+k\delta_2&\in J_2 \pmod{p_2^{3}},\\
&\setbox0\hbox{$\equiv$}\mathrel{\makebox[\wd0]{\vdots}}\\
\alpha_d+k\delta_d&\in J_d \pmod{p_d^{3}},
\end{align*}
where the sets $J_i$ consist of base-$p_i$ reversals of the numbers in the interval $[b_i,c_i)$ (which are $3$-digit long in base $p_i$). Let $L$ be the set of solutions of this system. Then $\cL'$ consists of all such sets $L$ as $(\alpha_i,\delta_i,b_i,c_i)_{i=1}^d$ ranges over all tuples in $(\mathbb{Z}/p_i^3 \mathbb{Z})^{4d}$.
\end{definition}
From the proof of \Cref{claim:L}, we have $L_B(\beta)\in \cL'$, for any $n$ and any good pair $(B,\beta)$. Also, note that $\abs{\cL'}\leq \gamma^{12}$. For each tuple $(\alpha_i,\delta_i,b_i,c_i)_{i=1}^d$, testing whether $k$ satisfies the equations can be done in $O(d)$ many arithmetic operations (by computing $d$ left-hand sides, reversing their digits, and seeing if the results are in appropriate intervals). Since any $L\in\cL'$ is $\gamma^3$-periodic, we only need to test $k$ satisfying $0\leq k<\gamma^3$. Thus, the set $\cL'$ can be computed in $O(d\gamma^{15})=d^{O(d)}$ many arithmetic operations.

\begin{definition}
We say that a subset $\bY'\subseteq [0,2\gamma^{15})$ is a \emph{representative} if there exist integers $\A',\D'$ and a set $L\in\cL'$ satisfying the following conditions:
\begin{itemize}
\item $\D'$ is even, and $\A'$ is divisible by $2\gamma^7$,
\item $0\leq \A'<\D'$,
\item $2\gamma^8-2< \D'\leq 8\gamma^{11}$,
\item $\bY'=(\A'+L\D')\cap [0,2\gamma^{15})$,
\item $|\bY'|\geq \gamma^4/16+1$.
\end{itemize}
\end{definition}

A representative is, roughly speaking, a sequence $\A(B)+L_B(\beta)\D(B)$ generated by the type of a good pair $(B,\beta)$, that is then scaled by $2\gamma^{11}/n$.

The input of the pre-processing step is just $d$, and the output is a set $X'\subseteq [0,2\gamma^{15})$ of size at most $900\gamma^3\log\gamma$ such that $X'+[0,\gamma^8)$ intersects all the representatives. 
The existence of $X'$ is guaranteed by the following two claims and union bound.
\begin{claim}\label{claim:repbound}
The number of representatives is at most $\gamma^{28}$.
\end{claim}
\begin{claim}\label{claim:hittingbound}
Let $x'$ be chosen uniformly from $[0,2\gamma^{15})$. For any fixed representative $\bY'$, the probability
of $[x',x'+\gamma^8)$ hitting $\bY'$ is
$\Pr\bigl[\,\bY\cap [x',x'+\gamma^8) \neq \emptyset \bigr]\geq 1/32\gamma^3$.
\end{claim}
\begin{proof}[Proof of~\Cref{claim:repbound}]
The numbers of possible $\A'$, $\D'$, and $L$ in the definition representatives are at most $8\gamma^{11}/2\gamma^7=4\gamma^4$, $8\gamma^{11}/2=4\gamma^{11}$, and $\gamma^{12}$, respectively. Thus, the number of representatives is at most $16\gamma^{27}\leq \gamma^{28}$.
\end{proof}
\begin{proof}[Proof of~\Cref{claim:hittingbound}] 
The proof is similar to the proof of \Cref{claim:probability}. Let $y\in \bY'$ be arbitrary. If $y\notin [0,\gamma^8)$, then $\Pr[y\in [x,x+\gamma^8)]=1/2\gamma^7$. Assume $\bY'$ is defined by $\A',\D',L$. Since
$\D'>2\gamma^8-2\geq \gamma^8$, the set $\bY'$ contains at most one element in the interval $[0,\gamma^8)$.
Hence \[\E\bigl[ \abs[\big]{\bY'\cap [x',x'+\gamma^8)} \bigr]\geq 1/32\gamma^3.\]
Since elements of $\bY'$ are at least $\D'$ apart, $\abs{\bY'\cap [x',x'+\gamma^8)}\in \{0,1\}$ for all $x'$.
Therefore, \[\Pr\bigl[\,\bY'\cap [x',x'+\gamma^8) \neq \emptyset \bigr]=\E\bigl[ \abs[\big]{\bY'\cap [x',x'+\gamma^8)} \bigr]\geq 1/32\gamma^3.\qedhere\]
\end{proof}

It follows from the second part of the claim that 
\begin{align*}
  \Pr\bigl[\,\bY'\cap \bigl(X'+[0,\gamma^8)\bigr)=\emptyset\bigr]\leq (1-1/32\gamma^3)^{900 \gamma^3 \log \gamma}<\gamma^{-28}.
\end{align*}
Thus, by the first part of the claim and union bound, there exists a choice of $X'$ such that $X'+[0,\gamma^8)$ intersects all the representatives. We can find such an $X'$ using the method of conditional expectations in time $\gamma^{O(1)}$.\smallskip

For any fixed $n$, which is divisible by $2\gamma^{11}$, we use $X'$ from the pre-processing step to find the set $P$ that meets every good box, similarly to the stage one of the construction in \Cref{thm:constr}. Namely, given $X'$ as above, let $X=\frac{n}{2\gamma^{11}}X'$. The desired set is then $P\eqdef r(X+[0,n/\gamma^3))$. Indeed, given a good pair $(B,\beta)$, let $\bY'$ be the representative 
corresponding to the triple $(\A',\D',L)$ where
$\A'=\A(B)\times 2\gamma^{11}/n$, $\D'=\D(B)\times 2\gamma^{11}/n$, and $L=\cL_B(\beta)$. From the pre-processing step, there exists some $x'\in X'$ such that $[x',x'+\gamma^8)\cap \bY'\neq \emptyset$. This implies $\A(B)+k\D(B)\in\frac{n}{2\gamma^{11}}x'+[0,n/2\gamma^3)$ for some $k\in L_B(\beta)$, and hence the set $P$ hits $\beta$.

Given such a set $P$ we can then proceed exactly as in the stage two in the proof of \Cref{thm:constr}. It is
completely deterministic. Naively, it takes $O(d\log n)$ steps
to compute each element of $P$ since computing the function $r_{p_i}$
requires $O(\log n)$ steps, for each $i$. That would make the total
number of operations in the algorithm $d^{O(d)}+O(dn\log n)$.
However, since the base-$p_i$ expansions of adjacent integers are
almost identical, it is possible to re-use the value of $r_{p_i}(m)$
when computing $r_{p_i}(m+1)$. This way one obtains an algorithm
with the total number of operations being $d^{O(d)}+O(dn)$.

\item Dispersion has also been studied on the torus. In this variant of the problem, 
the boxes are products of \emph{toroidal intervals}, which, in addition to the usual intervals $(a,b)$ for $a<b$,
include the sets of the form $(a,b)\eqdef (a,1]\cup [0,b)$ for $a>b$. Denote the $d$-dimensional
torus by $\Tor^d$, and let $m_d^{\Tor}$ be the corresponding dispersion function, i.e., the largest number
such that there is an empty box of volume $m_d^{\Tor}(n)$ among every $n$-point set on $\Tor^d$.
Ullrich \cite{ullrich} proved that $m_d^{\Tor}(n)\geq \min\{1,d/n\}$. This bound is trivially sharp for $d=1$,
and it was shown in \cite{breneis_hinrichs} that it is also sharp for $d=2$ and infinitely many~$n$.
In the opposite direction, the construction of Larcher, which was mentioned in the introduction, carries over verbatim to the
torus, and so $m_d(n)\leq 2^{7d+1}/n$. We can improve the base of exponent from $2^7$ to $e/2$.

\begin{proposition}
The toroidal dispersion satisfies $m_d^{\Tor}(n)\leq 32000(e/2)^dd^3\log d/n$, for all $n$ divisible by~$d$.
\end{proposition}
\begin{proof}
Let $P$ be the set obtained by invoking \Cref{thm:constr} with $n/d$ in place of~$n$.
Write $P+u$ to denote the shift of $P$ by vector $u$, where the `shift' is understood as a shift on~$\Tor^d$.
Set $v\eqdef\nobreak(1/d,1/d,...,1/d)\in \Tor^d$, and consider the shifts $P+rv$ for $r\in \{0,1,\dotsc,d-1\}$. 
We claim that the toroidal dispersion of $\bigcup_{r=0}^{d-1}(P+rv)$ is at most $32000(e/2)^dd^3\log d/n$. 
To prove this, it suffices, for every toroidal box $B_0$ of volume $32000(e/2)^dd^3\log d/n$, to find $r\in \{0,1,\dotsc,d-1\}$ such that the toroidal box
$B_r\eqdef B_0-rv$ contain a usual box of volume $8000d^3\log d/n$.
 
Write $\len(a,b)$ for the length of a toroidal interval $(a,b)$.
If $(a,b)$ is a toroidal interval, the largest usual interval contained in $(a,b)-x$ has length
$\len(a,b)$ if $x\notin (a,b)$ and $\max\bigl\{\len(a,x),\len(x,b)\bigr\}$ if $x\in (a,b)$.
For a toroidal interval $(a,b)$, let $f_{(a,b)}$ be the function given by 
\[
f_{(a,b)}(x)\eqdef
\begin{cases}
\log\left(\max\left\{\frac{\len(a,x)}{\len(a,b)},\frac{\len(x,b)}{\len(a,b)}\right\}\right)&\text{for }x\in(a,b),\\
 0&\text{for }x\notin (a,b).
\end{cases}
\]
If $B=\prod (a_i,b_i)$ is a toroidal box, the largest usual box contained in $B-(x_1,\dotsc,x_d)$ has volume 
\[
\vol(B)\exp\Bigl(\sum f_{(a_i,b_i)}(x_i)\Bigr).
\]

We shall estimate $\tfrac{1}{d}\sum_{r=0}^{d-1} f_{(a_i,b_i)}(r/d)$ by comparing it to the respective integral:
Since the function $f_{(a,b)}$ is unimodal with minimum at
$x=(a+b)/2$, the total variation of $f_{(a,b)}$ is $f_{(a,b)}(a)+f_{(a,b)}(b)-2f_{(a,b)}(\frac{a+b}{2})=2\log 2$.
Hence, 
\begin{equation}\label{eq:singleavg}
\frac{1}{d}\sum_{r=0}^{d-1}f_{(a,b)}(r/d)\geq \int_{0}^1f_{(a,b)}(x)\,dx-2\log2/d.
\end{equation} We can bound the integral in turn by
$
  \int f_{(a,b)}(x)\,dx=(\log 2-1)\len(a,b)\geq \log 2-1
$.
Summing~\eqref{eq:singleavg} over each of the $d$ coordinate directions, we then obtain
\[
\frac{1}{d}\sum_{r=0}^{d-1}\sum_{i=1}^d f_{(a_i,b_i)}(r/d)\geq d\log(2/e)-2\log 2.
\]
Hence, given any toroidal box $B_0$, there exists $r\in \{0,1,\dotsc,d-1\}$ such that the toroidal box $B_r=B_0-rv$ contains a usual box
of volume at least $\vol(B_0)\cdot \tfrac{1}{4}(2/e)^d$. In particular, if $\vol(B_0)\geq 32000(e/2)^dd^3\log d/n$,
then $B_r$ contains a usual box of volume $8000 d^3\log d/n$.
\end{proof}

It might be that the toroidal dispersion is indeed larger than the usual dispersion.
One evidence in that direction is that the VC dimension of boxes in the $[0,1]^d$
is $2d$ whereas the VC dimension of toroidal boxes is asymptotic to $d\log_2d$, as recently
showed by Gillibert, Lachmann and M\"ullner \cite{gillibert_lachmann_mullner}.

\item The first stage in the proof of \Cref{thm:constr} can be modified to yield a set $P$
such that the intersection $P\cap B$ with any dyadic box $B$ contains approximately
the same number of points. This can be used to give better constructions of sets in $\R^d$
without any large convex holes. The details are in \cite{almostnets} and \cite{convexholes}. 

\item We suspect that the smallest dispersion of an $n$-point set is asymptotic to
$\Theta(d\log d\cdot \frac{1}{n})$.

\end{itemize}

\bibliographystyle{plain}
\bibliography{emptybox}
\end{document}